\documentclass[preprint]{elsarticle} 

\usepackage{amssymb,amsmath,amsfonts, amscd, mathtools}
\usepackage{array}
\usepackage{enumerate}
\usepackage{graphicx}
\usepackage{color}
\usepackage{srcltx} 

\newtheorem{thm}{Theorem}
\newtheorem{fact}[thm]{Fact}
\newtheorem{conj}[thm]{Conjecture}
\newtheorem{cor}[thm]{Corollary}
\newtheorem{lem}[thm]{Lemma}
\newtheorem{defn}[thm]{Definition}
\newtheorem{prop}[thm]{Proposition}

\newtheorem{rem}[thm]{Remark}
\newcommand{\Khk}{K_h^k}

\let\eps=\varepsilon
\let\epsilon=\varepsilon

\newenvironment{pf}{\noindent\textbf{Proof.}}{\hfill$\Box$}
\newenvironment{pfcite}[1]{\noindent\textbf{Proof of~#1.}}{\hfill$\Box$}

\newcommand{\bdelta}{\hat{\delta}}
\newcommand{\textdef}{\textit}
\newcommand{\textindef}{\textbf}
\newcommand{\T}{\mathcal{T}}
\newcommand{\ess}{\mathcal{S}}

\begin{document}
\title{Asymptotic multipartite version of the Alon-Yuster theorem}
\author[isu]{Ryan R. Martin\fnref{fn1}}\ead{rymartin@iastate.edu}
\author[lse]{Jozef Skokan}\ead{j.skokan@lse.ac.uk}
\fntext[fn1]{This author's research partially supported by NSF grant DMS-0901008, NSA grant H98230-13-1-0226 and by an Iowa State University Faculty Professional Development grant.}
\address[isu]{Department of Mathematics, Iowa State University, Ames, Iowa 50011}
\address[lse]{Department of Mathematics, London School of Economics, London, WC2A 2AE, UK and Department of Mathematics, University of Illinois, 1409 W. Green Street, Urbana, IL 61801}

\begin{abstract}
	In this paper, we prove the asymptotic multipartite version of the Alon-Yuster theorem, which is a generalization of the Hajnal-Szemer\'edi theorem: If $k\geq 3$ is an integer, $H$ is a $k$-colorable graph and $\gamma>0$ is fixed, then, for every  sufficiently large $n$, where $|V(H)|$ divides $n$, and for every balanced $k$-partite graph $G$ on $kn$ vertices with each of its corresponding $\binom{k}{2}$ bipartite subgraphs having minimum degree at least $(k-1)n/k+\gamma n$, $G$ has a subgraph consisting of $kn/|V(H)|$ vertex-disjoint copies of $H$.

	The proof uses the Regularity method together with linear programming.
\end{abstract}

\maketitle

\section{Introduction}
\thispagestyle{empty}
\subsection{Motivation}
One of the celebrated results of extremal graph theory is the theorem of Hajnal and Szemer\'edi on tiling simple graphs with vertex-disjoint copies of a given complete graph $K_k$ on $k$ vertices. Let $G$ be a simple graph with vertex-set $V(G)$ and edge-set $E(G)$. We denote by $\deg_G(v)$, or simply $\deg(v)$, the degree of a vertex $v\in V(G)$ and we denote by $\delta(G)$ the minimum degree of the graph $G$. For a graph $H$ such that $|V(H)|$ divides $|V(G)|$, we say that $G$ has a \textdef{perfect $H$-tiling} (also a \textdef{perfect $H$-factor} or \textdef{perfect $H$-packing}) if there is a subgraph of $G$ that consists of $|V(G)|/|V(H)|$ vertex-disjoint copies of $H$.

The theorem of Hajnal and Szemer\'edi can be then stated in the following way:
\begin{thm}[Hajnal, Szemer\'edi~\cite{HSz}]\label{hsz}
	If $G$ is a graph on $n$ vertices, $k\mid n$, and $\delta(G)\geq (k-1)n/k$, then $G$ has a perfect $K_k$-tiling.
\end{thm}
The case of $k=3$ was first proven by Corr\'adi and Hajnal~\cite{CH} before the general case. The original proof in~\cite{HSz} was relatively long and intricate.  A shorter proof was provided later by Kierstead and Kostochka~\cite{KK}. Kierstead, Kostochka, Mydlarz and Szemer\'edi~\cite{KKMSz} improved this proof and gave a fast algorithm for finding $K_k$-tilings in $n$-vertex graphs with minimum degree at least $(k-1)n/k$.

The question of finding a minimum-degree condition for the existence of a perfect $H$-tiling in the case when $H$ is not a clique and $n$ obeys some divisibility conditions was first considered by Alon and Yuster~\cite{AY}:
\begin{thm}[Alon, Yuster~\cite{AY}]\label{ay}
	Let $H$ be an $h$-vertex graph with chromatic number $k$ and let $\gamma>0$. If $n$ is large enough, $h\mid n$ and $G$ is a graph on $n$ vertices with $\delta(G)\geq (k-1)n/k+\gamma n$, then $G$ has a perfect $H$-tiling.
\end{thm}
Koml\'os, S\'ark\"ozy and Szemer\'edi~\cite{KSSz01}
removed the $\gamma n$ term from the minimum degree condition and replaced it with a constant that depends only on~$H$.

K\"uhn and Osthus~\cite{KO} determined that $\left(1-1/\chi^*(H)\right)n+C$ was the necessary minimum degree to guarantee an $H$-tiling in an $n$-vertex graph for $n$ sufficiently large, and they also showed that this was best possible up to the additive constant.  The constant $C=C(H)$ depends only on $H$ and $\chi^*$ is an invariant related to the so-called critical chromatic number of $H$, which was introduced by Koml\'os~\cite{Kom}.


\subsection{Background}

In this paper, we consider the multipartite variant of Theorem~\ref{ay}. Before we can state the problem, we need a few definitions.

Given a graph $G$, the \textit{blow-up of $G$ by $m$}, denoted by $G(m)$, is the graph obtained by replacing each vertex $v\in V(G)$ with a set $U_v$ of $m$ vertices and replacing every edge $\{v_1,v_2\}\in E(G)$ with the complete bipartite graph $K_{m,m}$ on vertex sets $U_{v_1}$ and $U_{v_2}$. 

   A $k$-partite graph $G=(V_1,\ldots,V_k; E)$ is \textit{balanced} if $|V_1|=\cdots =|V_k|$. The \textit{natural bipartite subgraphs} of $G$ are those induced by the pairs $(V_i, V_j)$, and which we denote by $G[V_i, V_j]$. For a $k$-partite graph $G=(V_1,\ldots,V_k;E)$, we define the \textit{minimum bipartite degree}, $\bdelta_k(G)$, to be the smallest minimum degree among all of the natural bipartite subgraphs of $G$, that is,
$$	\bdelta_k(G)=\min_{1\leq i < j \leq k}\delta(G[V_i, V_j]) .	$$

Now we can state the conjecture that inspired this work, a slightly weaker version of which appeared in~\cite{CsM}.
\begin{conj}\label{conj:main}
	Fix an integer $k\geq 3$. If $G$ is a balanced $k$-partite graph on $kn$ vertices such that $\bdelta_k(G)\geq (k-1)n/k$, then either $G$ has a perfect $K_k$-tiling or both $k$ and $n/k$ are odd integers and $G$ is isomorphic to the fixed graph $\Gamma_{k,n}$.
\end{conj}
The exceptional graphs $\Gamma_{k,n}$, where $n$ is an integer divisible by $k$, are due to Catlin~\cite{C} who called them ``type 2 graphs''. The graph $\Gamma_{k,k}$ has vertex set $\left\{h_{ij}:i,j\in\{1,\ldots,k\}\right\}$ and $h_{ij}$ is adjacent to $h_{i'j'}$ if $i\neq i'$ and either $j=j'\in \{k-1, k\}$ or $j\neq j'$ and at least one of~$j, j'$ is in $\{1,\ldots,k-2\}$.  For $n$ divisible by $k$,  the graph $\Gamma_{k,n}$ is the blow-up graph $\Gamma_{k,k}(n/k)$.

We notice that if $G$ satisfies the minimum bipartite degree condition in Conjecture~\ref{conj:main}, then its minimum degree $\delta(G)$ can still be as small as $(k-1)\left(\frac{k-1}{k}\right)n=\left(\frac{k-1}{k}\right)^2(kn)$, which is not enough to apply Theorem~\ref{hsz} directly.

The case of $k=2$ of Conjecture~\ref{conj:main} is an immediate corollary of the classical matching theorem due to K\"onig~\cite{Kon} and Hall~\cite{Hal}. Fischer~\cite{Fis} observed that if $G$ is a balanced $k$-partite graph on $kn$ vertices with $\bdelta_k(G)\geq\left(1-1/2(k-1)\right)n$, then $G$ has a perfect $K_k$-tiling. 

Some partial results were obtained, for $k=3$, by Johansson~\cite{Joh} and, for $k=3,4$, by Fischer~\cite{Fis}. The case of $k=3$ was settled for $n$ sufficiently large by Magyar and the first author~\cite{MM}, and the case of $k=4$ was settled for $n$ sufficiently large by Szemer\'edi and the first author~\cite{MSz}. The results in~\cite{MM,MSz} each have as a key lemma a variation of the results of Fischer.  However, it seems that such techniques are impossible for $k\geq 5$. An interesting result toward proving Conjecture~\ref{conj:main} for general $k$ is due to Csaba and Mydlarz~\cite{CsM} who proved that if $G$ is a balanced $k$-partite graph on $kn$ vertices, $\bdelta_k(G)\geq\tfrac{q_k}{q_k+1} n$ and $n$ is large enough, then $G$ has a perfect $K_k$-tiling. Here, $q_k:=k-\tfrac{3}{2}+\tfrac{1}{2}\sum\limits_{i=1}^k\tfrac{1}{i}=k+O(\log k)$.

Recently, Keevash and  Mycroft~\cite{KM} proved that, for any $\gamma>0$, if $n$ is large enough, then $\bdelta_k(G)\geq(k-1)n/k+\gamma n$ guarantees a perfect $K_k$-tiling in a balanced $k$-partite graph $G$ on $kn$ vertices. Their result is a consequence of a more general theorem on hypergraph matching, the proof of which uses the hypergraph regularity method and a hypergraph version of the Blow-up Lemma. Very shortly thereafter, Lo and Markstr\"om \cite{LM} proved the same result using methods from linear programming and the so-called ``absorbing method''. This effort culminated in \cite{KM2}, in which Keevash and  Mycroft proved Conjecture~\ref{conj:main}.

In this paper, we are interested in more general problem of tiling $k$-partite balanced graphs by a fixed $k$-colorable graph $H$. More precisely, if $H$ is a $k$-colorable graph and $n$ obeys certain natural divisibility conditions, we look for a condition on $\bdelta_k(G)$ to ensure that every balanced $k$-partite graph $G$ on $kn$ vertices satisfying this condition has a perfect $H$-tiling.

Zhao~\cite{Z} found that the minimum degree required to perfectly tile a balanced bipartite graph on $2n$ vertices with copies of $K_{h,h}$ ($h$ divides $n$) is $n/2+C(h)$, where $C(h)$ differs sharply as to whether $n/h$ is odd or even.  Zhao and the first author~\cite{MZ1,MZ2} showed similar results for tiling with $K_{h,h,h}$. Hladk\'y and Schacht~\cite{HS} and then Czygrinow and DeBiasio~\cite{CDeB} improved the results of~\cite{Z} by finding the minimum degree for copies of $K_{s,t}$, where $s+t$ divides $n$. Bush and Zhao~\cite{BZ} proved a K\"uhn-Osthus-type result by finding the asymptotically best-possible minimum degree condition in a balanced bipartite graph on $2n$ vertices in order to ensure its perfect $H$-tiling, for any bipartite $H$.  All results are for $n$ sufficiently large.

\subsection{Main Result}
We prove a multipartite version of the Alon-Yuster theorem (Theorem~\ref{ay}).  Let $\Khk$ denote a $k$-partite graph with $h$ vertices in each partite set. For example, the complete bipartite graph $K_{h,h}$ would be denoted $K_h^2$.  Since the partite sets can be rotated, it is easy to see that any $k$-chromatic graph $H$ of order $h$ perfectly tiles the graph $\Khk$. Hence, the following theorem gives a sufficient condition for a perfect $H$-tiling.
\begin{thm}\label{thm:main}
	Fix an integer $k\geq 2$, an integer $h\geq 1$ and $\gamma\in (0,1)$. If $n$ is sufficiently large, divisible by $h$, and $G$ is a balanced $k$-partite graph on $kn$ vertices with $\bdelta_k(G)\geq\left(\frac{k-1}{k}+\gamma\right)n$, then $G$ has a perfect $\Khk$-tiling.
\end{thm}
Our proof relies on the regularity method for graphs and linear programming and it differs from approaches in \cite{KM,LM}.

\subsection{Structure of the Paper}
In Section~\ref{sec:LP}, we prove a fractional version of the multipartite Hajnal-Szemer\'edi theorem.  This is the main tool in proving Theorem~\ref{thm:main}.  Section~\ref{sec:mainproof} is the main proof and Section~\ref{sec:facts} gives the proofs of the supporting lemmas.  We finish with Section~\ref{sec:conc}, which has some concluding remarks.

\section{Linear Programming}
\label{sec:LP}
\newcommand*\boxSizeOfmax[1]{\makebox[\widthof{$\max$}][c]{#1}}
\newcommand*\boxSizeOfmin[1]{\makebox[\widthof{$\min$}][c]{#1}}

In this section, we shall prove a fractional version of Conjecture~\ref{conj:main}.
\begin{defn}
	For any graph $G$, let $\T_k(G)$ denote the set of all copies of $K_k$ in $G$.  The \textindef{fractional $K_k$-tiling number} $\tau_k^*(G)$ is defined as:
	\begin{align}\label{eq:maxLP}
   		\tau_k^*(G) 	=	& \max 				& \sum_{T\in\T_k(G)}w(T) \\
   						& \boxSizeOfmax{\rm s.t.} 	& \sum_{\scriptsize\begin{array}{c}T\in\T_k(G) \\ V(T)\ni v \end{array}}w(T) 
						& \leq 1, 	& \quad\forall v\in V(G), \nonumber \\
   						& 					& w(T) 
						& \geq 0, 	& \quad\forall T\in\T_k(G) .  \nonumber
	\end{align}
\end{defn}
From the Duality Theorem of linear programming (see~\cite[Section 7.4]{Sch}), we obtain that
\begin{align}\label{eq:minLP}
	\tau_k^*(G)			 =	& \min 				& \sum_{v\in V(G)}x(v) \\
                  				& \boxSizeOfmin{s.t.} 		& \sum_{v\in V(T)}x(v) 
						& \geq 1, 				& \quad \forall T\in\T_k(G), \nonumber \\
                            			& 					& x(v) 
						& \geq 0, 				& \quad \forall v\in V(G). \nonumber
\end{align}
 
   Let $w^*$ be a function that achieves an optimal solution to (\ref{eq:maxLP}). If there exists a vertex $v\in V(G)$ such that $\sum_{T\in\T_k(G), V(T)\ni v}w^*(T)<1$, then we call $v$ a \textdef{slack vertex} or just say that \textdef{$v$ is slack}.  Similarly, if $x^*$ is a function that achieves an optimal solution to (\ref{eq:minLP}) and there exists a $T\in\T_k(G)$ such that $\sum_{v\in V(T)}x^*(v)>1$, then we say that \textdef{$T$ is slack}.
 
\begin{rem}\label{rem:rational}
	Consider an optimal solution to (\ref{eq:maxLP}), call it $w^*$.  We may assume that $w^*(T)$ is rational for each $T\in\T_k(G)$.  To see this, observe that the set of feasible solutions is a polyhedron for which each vertex is the solution  to a system of equations that result from setting a subset of the constraints of the program (\ref{eq:maxLP}) to equality. (For more details, see~\cite[Theorem 18.1]{Chv}.)  Since the objective function achieves its maximum at such a vertex (See~\cite[Section 3.2]{Gass}.) we may choose an optimal solution $w^*(T)$ with rational entries.
\end{rem}

Now we can state and prove a fractional version of the multipartite Hajnal-Szemer\'edi Theorem.
\begin{thm}\label{thm:frac}
	Let $k\geq 2$. If $G$ is a balanced $k$-partite graph on $kn$ vertices such that $\bdelta_k(G)\geq(k-1)n/k$, then $\tau_k^*(G)=n$.
\end{thm}

\begin{pf}
	Setting $x(v)=1/k$ for all vertices $v\in V(G)$ gives a feasible solution $x$ to (\ref{eq:minLP}), and so $\tau_k^*(G)\leq \sum\limits_{v\in V(G)}x(v)=n$.  We establish that $\tau_k^*(G)\geq n$ by induction on $k$. \\

	\noindent\textbf{Base Case. $k=2$.}
	This case follows from the fact that Hall's matching condition implies that a balanced bipartite graph on $2n$ vertices with minimum degree at least $n/2$ has a perfect matching.  Setting $w(e)$ equal to $1$ if edge $e$ is in the matching and equal to $0$ otherwise, gives a feasible solution to (\ref{eq:maxLP}), thus establishing that $\tau_2^*(G)\geq n$. \\

	\noindent\textbf{Induction step. $k\geq 3$.} Now we assume $k\geq 3$ and suppose, for any balanced $(k-1)$-partite graph $G'$ on a total of $(k-1)n'$ vertices with $\bdelta_{k-1}(G')\geq\frac{k-2}{k-1}n'$, that $\tau_{k-1}^*(G')\geq n'$.

	Let $w^*$ be an optimal solution to (\ref{eq:maxLP}).  Let $x^*$ be an optimal solution corresponding to (\ref{eq:minLP}) such that $x^*(z)=0$ whenever vertex $z$ is slack. This is guaranteed by the Complementary Slackness Theorem~\cite[Section 7.9]{Sch}.  Denote by $\ess$ the set of slack vertices, and, for $i\in [k]$, set $\ess_i=\ess\cap V_i$.   If some $\ess_i=\emptyset$, then $V_i$ having no slack vertices gives that $\sum_{T\ni v}w^*(T)=1$ for each $v\in V_i$. Since each $T\in\T_k(G)$ has exactly one vertex in $V_k(G)$, then $\tau_k^*(G)=n$. Hence, we may assume that every $\ess_i$ is non-empty.

	Denote $[k]:=\{1,\ldots,k\}$. For every $i\in [k]$, fix some $z_i\in\ess_i$, choose exactly $n':=\left\lceil\frac{k-1}{k}n\right\rceil$ neighbors of $z_i$ in each $V_j$, $j\in [k]-\{i\}$, and denote by $G_i$ the subgraph of $G$ induced on these $(k-1)n'$ neighbors.

	Observe that the set of weights $\{x^*(v) : v\in V(G_i)\}$ must be a feasible solution to the minimization problem (\ref{eq:minLP})  defined by the $(k-1)$-partite graph $G_i$. This is because every copy of $K_{k-1}$ in $G_i$ extends to a copy of $K_k$ in $G$ containing the vertex $z_i$ and the sum of the weights of the vertices on that $K_{k-1}$ must be at least $1$ because $x^*(z_i)=0$. Hence, we have that $\sum\limits_{v\in V(G_i)}x^*(v)\geq\tau^*_{k-1}(G_i)$.

	Each vertex of $G_i$ has at most $n-n'$ neighbors outside of $V(G_i)$ in each of its classes.  Thus,
	$$	\bdelta_{k-1}(G_i)\geq n'-(n-n')=\frac{k-2}{k-1}n'+\left(\frac{k}{k-1}n'-n\right)\geq \frac{k-2}{k-1}n' . 	$$
  	So, for every $i$, we may apply the inductive hypothesis to $G_i$ and conclude that $\tau^*_{k-1}(G_i)=n'$.

	Combining the previous two observations with the fact that each vertex $v$ is in at most $k-1$ of the subgraphs $G_i$, we get
   	$$ 	(k-1)\tau_k^*(G) = (k-1)\sum_{v\in V(G)}x^*(v) \geq \sum_{i=1}^k\sum_{v\in V(G_i)}x^*(v) \geq \sum_{i=1}^k\tau^*_{k-1}(G_i) = kn' . 	$$

	So, $\tau_k^*(G)\geq\frac{k}{k-1}n'=\frac{k}{k-1}\left\lceil\frac{k-1}{k}n\right\rceil\geq n$.  This concludes the proof of Theorem~\ref{thm:frac}.
\end{pf}~\\

\section{Proof of Theorem~\ref{thm:main}}
\label{sec:mainproof}
First, we will have a sequence of constants and the notation $a\gg b$ means that the constant $b$ is sufficiently small compared to $a$. We fix $k\geq 2$ and $h\geq 1$ and let
\begin{equation}\label{eq:epsprime}
	\min\{k^{-1},h^{-1},\gamma\}\gg d\gg\epsilon'\gg\zeta\gg n^{-1} ,
\end{equation}
We have an additional parameter $\epsilon$ and specify that $\epsilon=(\epsilon')^5/16$.

\subsection{Applying the Regularity Lemma}

We are going to use a variant of Szemer\'edi's Regularity Lemma.  Before we can state it, we need a few basic definitions. If $G$ is a graph with $S\subset V(G)$ and $x\in V(G)$, then $\deg_G(x,S)$ (or $\deg(x,S)$ if $G$ is understood) denotes $|N(x)\cap S|$.

For disjoint vertex sets $A$ and $B$ in some graph, let $e(A,B)$ denote the number of edges with one endpoint in $A$ and the other in $B$. Further, let the \textdef{density} of the pair $(A,B)$ be $d(A,B)=e(A,B)/|A||B|$. The pair $(A,B)$ is \textdef{$\epsilon$-regular} if $X\subseteq A$, $Y\subseteq B$, $|X|\geq\epsilon |A|$ and $|Y|\geq\epsilon|B|$ imply $|d(X,Y)-d(A,B)|\leq\epsilon$.

We say that a pair $(A,B)$ is \textdef{$(\epsilon,\delta)$-super-regular} if it is $\epsilon$-regular and $\deg(a,B)\geq\delta |B|$ for all $a\in A$ and $\deg(b,A)\geq\delta |A|$ for all $b\in B$.

The degree form of Szemer\'edi's Regularity Lemma (see, for instance,~\cite{KS}) is sufficient here, modified for the multipartite setting.

\begin{thm}\label{thm:SzemRegLem}
	For every integer $k\geq 2$ and every $\epsilon>0$, there is an $M=M(k,\epsilon)$ such that if $G=(V_1,\ldots,V_k;E)$ is a balanced $k$-partite graph on $kn$ vertices and $d\in[0,1]$ is any real number, then there is an integer $\ell$, a subgraph $G'=(V_1,\ldots,V_k;E')$ and, for $i=1,\ldots,k$, partitions of $V_i$ into clusters $V_i^{(0)},V_i^{(1)},\ldots,V_i^{(\ell)}$ with the following properties:
	\begin{enumerate}[(P1)]
		\item $\lceil\epsilon^{-1}\rceil\leq\ell\leq M$,
		\item $|V_i^{(0)}|\leq \epsilon n$ for $i\in [\ell]$,
		\item $|V_i^{(j)}|=L\leq\epsilon n$ for $i\in [k]$ and $j\in [\ell]$,
		\item $\deg_{G'}(v,V_{i'})>\deg_G(v,V_{i'})-(d+\epsilon)n$ for all $v\in V_i$, $i\neq i'$, and
		\item all pairs $(V_i^{(j)},V_{i'}^{(j')})$, $i,i'\in [k]$, $i\neq i'$, $j,j'\in [\ell]$, are $\epsilon$-regular in $G'$, each with density either $0$ or exceeding $d$. \label{it:P5}
	\end{enumerate}
\end{thm}
We omit the proof of Theorem~\ref{thm:SzemRegLem}, which follows from the proof given in~\cite{Szem}.

Given a balanced $k$-partite graph $G$ on $kn$ vertices with $\bdelta_k(G)\geq\left(\frac{k-1}{k}+\gamma\right)n$, and given $d$ and $\epsilon$, we construct the \textdef{reduced graph} $G_r$ on $k\ell$ vertices corresponding to the clusters $V_i^{(j)}$, $1\leq i\leq k$, $1\leq j\leq\ell$, obtained from Theorem~\ref{thm:SzemRegLem}.  Each edge of $G_r$ corresponds to an $\epsilon$-regular pair with density at least $d$ in $G'$. Observe that $G_r$ is $k$-partite and balanced. Lemma~\ref{lem:mindeg} shows that $G_r$ has a similar minimum-degree condition to that of $G$.
\begin{lem}\label{lem:mindeg}
	Let $G$ be a balanced $k$-partite graph $G$ on $kn$ vertices with $\bdelta_k(G)\geq\left(\frac{k-1}{k}+\gamma\right)n$. Then, for the reduced graph $G_r$ defined as above, we have $\bdelta_k(G_r)\geq\left(\frac{k-1}{k}+\gamma-((k+2)\epsilon+d)\right)\ell$. Furthermore, if $(k+2)\epsilon+d\leq\gamma/2$, then
	$$ 	\bdelta_k(G_r)\geq\left(\frac{k-1}{k}+\gamma/2\right)\ell . 	$$
\end{lem}

The proof of Lemma~\ref{lem:mindeg} is immediate (see~\cite{CsM}).~\\

\subsection{Partitioning the clusters}
\label{sec:partitioning}
We first apply the fractional version of the $k$-partite Hajnal-Szemer\'edi Theorem (Theorem~\ref{thm:frac}) to $G_r$ and obtain that the value of $\tau_k^*(G_r)$ is equal to $\ell$.  Consider a corresponding optimal solution $w^*$ to the linear program (\ref{eq:maxLP}) as it is applied to $G_r$. By Remark~\ref{rem:rational}, we may fix a corresponding solution $w^*(T)$ that is rational for every $T\in\T_k(G_r)$. We will call this $w^*$ a \textdef{rational-entry solution for $G_r$} and denote by $D(G_r)$ the common denominator of all of the entries of~$w^*$. 

Since the linear program (\ref{eq:maxLP}) depends only on $G_r$ and the number of such reduced graphs is only dependent on $M(k,\epsilon)$,
the number of possible linear programs is only dependent only on $k$ and $\epsilon$. For each possible linear program we fix one rational-entry solution.

Therefore, the least common multiple of all of the common denominators $D(G_r)$ for these reduced graphs is a function only of $k$ and $\epsilon$.  Call it $D=D(k,\epsilon)$.  In sum, $D$ has the property that for every reduced graph $G_r$, there is a rational-entry solution $w^*$ of the linear program (\ref{eq:maxLP}) such that $D\cdot w^*(T)$ is an integer for every $T\in\T_k(G_r)$.

The next step is to partition, uniformly at random, each set $V_i^{(j)}$ into $D$ parts of size $h\lfloor L/(Dh)\rfloor$ as well as a single (possibly empty) set of size $L-Dh\lfloor L/(Dh)\rfloor<Dh$. The vertices of the latter set of less than $Dh$ vertices will be added to the corresponding leftover set, $V_i^{(0)}$.  The resulting leftover set $\tilde{V}_i^{(0)}$ has size less than $\epsilon n+Dh\ell<2\epsilon n$.

Thus, for $L'=h\lfloor L/(Dh)\rfloor$, we obtain $k(D\ell)$ clusters $\tilde{V}_i^{(j)}$, $i\in [k]$, $j\in [D\ell]$, such that each of them has size exactly $L'$. This new partition has the following properties:
\begin{enumerate}[(P1$'$)]
	\item $\ell'=D\ell$, \label{it:P1prime}
	\item $|\tilde{V}_i^{(0)}|\leq 2\epsilon n$ for $i\in [k]$,\label{it:P2prime}
	\item $|\tilde{V}_i^{(j)}|=L'=h\lfloor L/(Dh)\rfloor$ for $i\in [k]$ and $j\in [\ell']$, \label{it:P3prime}
	\item $\deg_{G'}(v,V_{i'})>\deg_G(v,V_{i'})-(d+\epsilon)n$ for all $i,i'\in [k]$, $i\neq i'$, $v\in V_i$ and \label{it:P4prime}
\end{enumerate}

Now we prove that a property similar to property (P\ref{it:P5}) holds.
\begin{enumerate}[(P1$'$)]\setcounter{enumi}{4}
	\item all pairs $(\tilde{V}_i^{(j)},\tilde{V}_{i'}^{(j')})$, $i,i'\in [k]$, $i\neq i'$, $j,j'\in [\ell']$ are $\epsilon'$-regular in $G'$, each with density either $0$ or exceeding $d':=d-\epsilon$.\label{it:P5prime}
\end{enumerate}
Recall from (\ref{eq:epsprime}) that $\epsilon=(\epsilon')^5/16$ and, consequently, $\epsilon'=(16\epsilon)^{1/5}$.

The upcoming Lemma~\ref{lem:randslice}, a slight modification of a similar lemma by Csaba and Mydlarz~\cite[Lemma 14]{CsM}, implies that, in fact, (P\ref{it:P5prime}$'$) holds with probability going to 1 as $n\rightarrow\infty$. The proof follows easily from theirs and so we omit it.

\begin{lem}[Random Slicing Lemma]\label{lem:randslice}
	Let $0<d<1$, $0<\epsilon<\min\{d/4,(1-d)/4,1/9\}$ and $D$ be a positive integer.  There exists a $C=C(\epsilon,D)>0$ such that the following holds: Let $(X,Y)$ be an $\epsilon$-regular pair of density $d$ with $|X|=|Y|=DL'$. If $X$ and $Y$ are randomly partitioned into sets $A_1,\ldots,A_D$, and $B_1,\ldots,B_D$, respectively, each of size $L'$, then, with probability at least $1-\exp\{-C\cdot DL'\}$, all pairs $(A_i,B_j)$ are $(16\epsilon)^{1/5}$-regular with density at least $d-\epsilon$.
\end{lem}

Using Lemma~\ref{lem:randslice}, the property (P\ref{it:P5prime}$'$) holds with probability at least $1-\binom{k}{2}\ell^2\exp\{-CDL'\}=1-\binom{k}{2}\ell^2\exp\{-O(L)\}$. Since $\ell\leq M=M(k,\epsilon)$ and $L\geq n(1-\epsilon)/M$, then for every sufficiently large~$n$, a partition satisfying (P\ref{it:P1prime}$'$)-(P\ref{it:P5prime}$'$) exists (with high probability). We fix a partition that satisfies (P\ref{it:P1prime}$'$)-(P\ref{it:P5prime}$'$). The sets $\tilde{V}_j^{(j)}$ are called \textdef{sub-clusters}.

To understand this new partition, we define its reduced graph $G_r'$ with vertex set $\bigcup_{i=1}^k \{u_i^{(1)},\ldots,u_i^{(\ell')}\}$. The vertex $u_i^{(j)}$ corresponds to the cluster $\tilde{V}_i^{(j)}$. The vertices $u_i^{(j)}$ and $u_{i'}^{(j')}$ are adjacent in $G_r'$ if and only if the pair $(\tilde{V}_i^{(j)},\tilde{V}_{i'}^{(j')})$ is $\epsilon'$-regular with density at least $d'$.  The graph $G_r'$ clearly has the following properties:
\begin{itemize}
	\item $G_r'$ is $k$-partite and balanced on $k\ell'$ vertices. We denote its partite sets  $U'_i=\{u_i^{(1)},\ldots,u_i^{(\ell')}\}$, $i\in [k]$.
	\item $\bdelta_k(G_r')\geq\left(\frac{k-1}{k}+\gamma/2\right)\ell'$.
\end{itemize}

The usefulness of $G_r'$ is that it has a $K_k$-tiling, which is derived from the fractional $K_k$-tiling of $G_r$:

\begin{fact}\label{fact:grtiling}
	The reduced graph $G_r'$ has a perfect $K_k$-tiling.
\end{fact}

\begin{pfcite}{Fact~\ref{fact:grtiling}}
	Observe first that, by (P\ref{it:P5}) and (P\ref{it:P5prime}$'$), $G_r'$ is simply the blow-up graph $G_r(D)$. Let $w^*$ be the previously-chosen rational-valued solution to the linear program (\ref{eq:maxLP}) as applied to $G_r$.

	Consider some $T\in\T_k(G_r)$ with vertices $\{v_1,\ldots,v_k\}$. Observe that, by the definition of $D$, $D w^*(T)$ is an integer.  Then, we take $D w^*(T)$ of the vertices from $U_{v_1}$, $D w^*(T)$ of the vertices from $U_{v_2}$ and so on until taking $D w^*(T)$ of the vertices from $U_{v_k}$.  This  selection produces $D w^*(T)$ vertex-disjoint copies of $K_k$ in $G_r'$.

	By the constraint inequalities in (\ref{eq:maxLP}), the total number of vertices used from $U_v$ is
	$$ 	\sum_{T\in\T_k(G_r), V(T)\ni v}D w^*(T)\leq D=|U_v| , 	$$
  	hence the process never fails.  The total number of vertex-disjoint $K_k$-s that are created in this way is $\sum_{T\in\T_k(G_r)}D w^*(T)=D \ell=\ell'$. This uses each of the $k\ell'$ vertices of $G_r'$.
\end{pfcite}~\\

Since $G_r'$ has a perfect tiling, we may re-index its vertices so that vertices of $G_r'$ (the vertices of $G_r'$ correspond to the sub-clusters of $G$) with the same upper-index are in the same copy of the tiling from Fact~\ref{fact:grtiling}.  More precisely,
\begin{itemize}
	\item for $j=1,\ldots,\ell'$, the $k$-tuple $(u_1^{(j)},\ldots,u_k^{(j)})$ forms a $K_k$ in $G_r'$. We refer to the $k$-tuples $(\tilde{V}_1^{(j)},\ldots,\tilde{V}_k^{(j)})$ as \textdef{columns}.\footnote{We visualize the vertex sets $V_i$ as being horizontal, like rows in a matrix, so it is natural to think of these $k$-tuples as columns.}
\end{itemize}

\subsection{Making the cliques super-regular}

In preparation for using the Blow-up Lemma (Lemma~\ref{lem:blow-up} below), we need to make each $k$-tuple $(\tilde{V}_1^{(j)},\ldots,\tilde{V}_k^{(j)})$, $j\in [\ell']$, pairwise super-regular by placing some vertices from the corresponding sub-clusters into the respective leftover set. This is easy to do by a simple fact which is proven in Section~\ref{sec:facts}:
\begin{fact}\label{fact:super-slicing}
	Let $\epsilon'>0$ and $\epsilon'<d'/(2k+2)$. Let $(A_1,\ldots,A_k)$ be a $k$-tuple that is pairwise $\epsilon'$-regular of density at least $d'$ with $|A_1|=\cdots=|A_k|=L'$.  There exist subsets $A_i'\subset A_i$ for $i\in [k]$ such that $|A_i|=h\lceil (1-(k-1)\epsilon')L'/h\rceil$ and each pair of $(A_1',\ldots,A_k')$ is $(2\epsilon',d'-k\epsilon')$-super-regular (with density at least $d'-\epsilon'$).
\end{fact}

Fact~\ref{fact:super-slicing} follows from well-known properties of regular pairs. We apply it to each $k$-tuple $(\tilde{V}_1^{(j)},\ldots,\tilde{V}_k^{(j)})$, $j\in [\ell']$.  We do not rename the sets $\tilde{V}_i^{(j)}$ since they only shrink in magnitude only by $(k-1)\epsilon'L'$. Consequently,
\begin{itemize}
	\item the leftover sets $\tilde{V}_i^{(0)}$, $1\leq i\leq k$, are of size at most $2\epsilon n+(k-1)\epsilon'L'\ell'<k\epsilon' n$,
	\item each pair $(\tilde{V}_i^{(j)},\tilde{V}_{i'}^{(j)})$, $i\neq i'$, is $(2\epsilon',d'/2)$-super-regular, and
	\item each pair $(\tilde{V}_i^{(j)},\tilde{V}_{i'}^{(j')})$ is $2\epsilon'$-regular with density either 0 or at least $d'-\epsilon'$, regardless of whether or not $j=j'$.
\end{itemize}

If we use the Blow-up Lemma (Lemma~\ref{lem:blow-up}) at this point, we would obtain a $\Khk$-tiling that covers every vertex of $G$ except those in the leftover sets. The remainder of the proof is to establish that we can, in fact, ensure that the leftover vertices can be absorbed by the sub-clusters and we can obtain a $\Khk$-tiling that covers all the vertices of $G$.

\subsection{Preparing for absorption}
\label{sec:prepabs}

In order to absorb the vertices from the leftover sets, we need to prepare some copies of $\Khk$ throughout $G$ that may be included in the final $\Khk$-tiling.  Their purpose is to ensure that, after inserting vertices from the leftover sets to the sub-clusters, the number of vertices in each of the sub-clusters can be balanced so that the Blow-up Lemma (Lemma~\ref{lem:blow-up}) can be used.  These copies of $\Khk$ will be specially designated and colored either red or blue according to their role.

The Reachability Lemma (Lemma~\ref{lem:reach}) is how we transfer the imbalance of the sizes of one column to the first column.

\begin{figure}[ht]
	\begin{center}
		\includegraphics[height=2in]{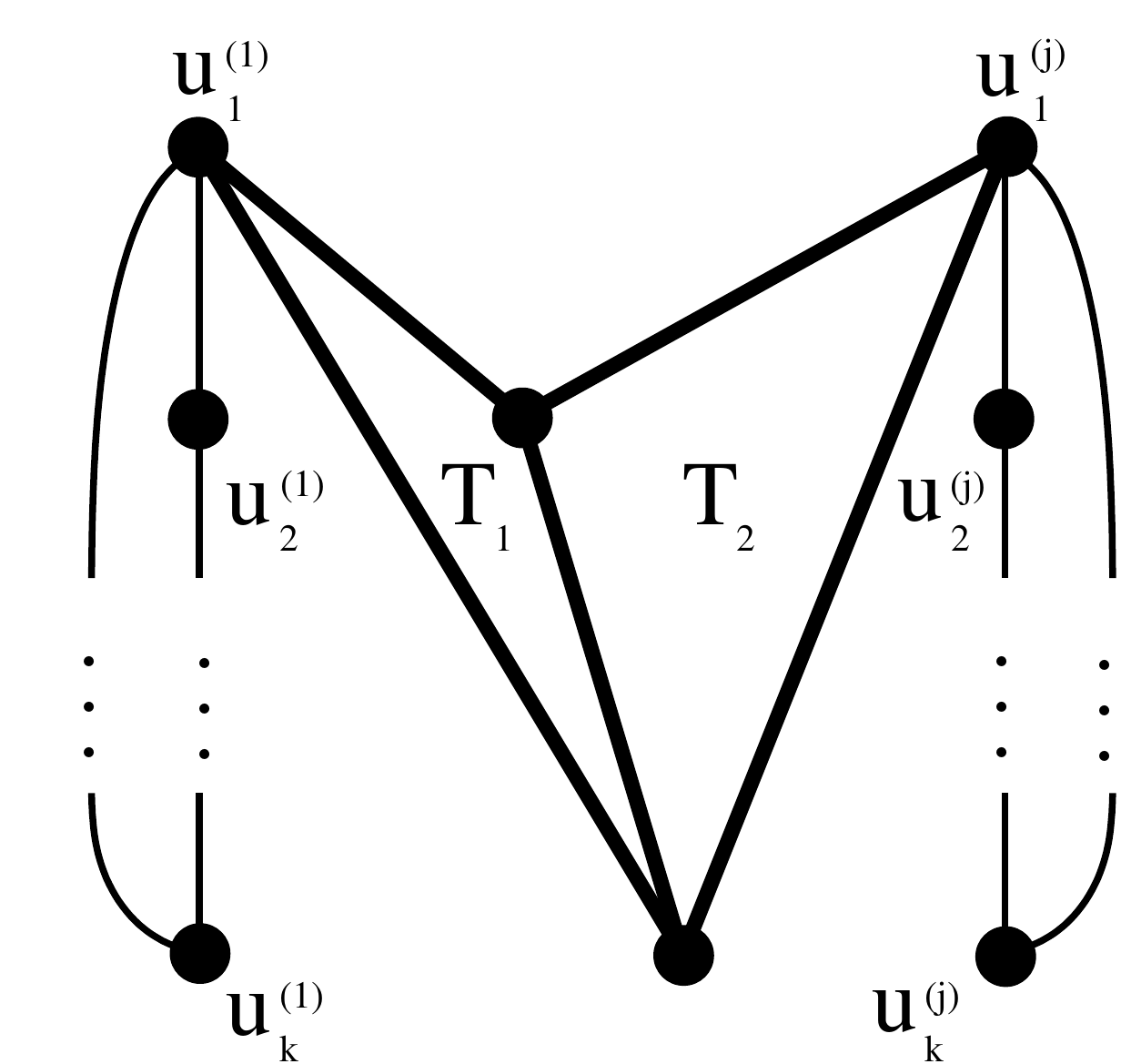}
		\caption{Diagram for $T_1$ and $T_2$ formed in reaching $u_1^{(1)}$ from $u_1^{(j)}$.}
	\end{center}
	\label{fig1}
\end{figure}

\begin{lem}[Reachability Lemma]\label{lem:reach}
	Let $G_r'$ be a balanced $k$-partite graph with partite sets $U_i'=\{u_i^{(j)} : j\in [\ell']\}$, $i\in [k]$. Let $\bdelta_k(G_r')\geq\frac{k-1}{k}\ell'+2$. Then, for each $i\in [k]$ and $j\in\{2,\ldots,\ell'\}$, there is a pair $(T_1,T_2)$ of copies of $K_k$  such that their symmetric difference is $\{u_i^{(1)},u_i^{(j)}\}$ and $T_1$ and $T_2$ contain no additional vertices from $\{u_1^{(1)},\ldots,u_k^{(1)},u_1^{(j)},\ldots,u_k^{(j)}\}$. See Figure~\ref{fig1}.
\end{lem}

\begin{pfcite}{Lemma~\ref{lem:reach}}
	Without loss of generality, it suffices to prove the lemma for $i=1$ and $j=\ell'$.  The vertices $u_1^{(1)}$ and $u_1^{(\ell')}$ have at least $\ell'-2(\ell'-\bdelta_k(G_r'))\geq\ell'-2\left(\ell'-\frac{k-1}{k}\ell'-2\right)=\left(\frac{k-2}{k}\right)\ell'+4$ common neighbors in each of $U_2',\ldots,U_k'$. Hence, one can choose a sequence $w_2,\ldots,w_k$ of vertices so that, for $i=2,\ldots,k$, $w_i$ is in $U_i'-\{u_i^{(1)},u_i^{(\ell')}\}$ and is a common neighbor of $u_1^{(1)},u_1^{(\ell')},w_2,\ldots,w_{i-1}$. Note that at each stage, the number of available choices for $w_i$ is at least
	$$	 \left(\frac{k-2}{k}\ell'+4\right)-(i-2)\left(\ell'-\frac{k-1}{k}\ell'-2\right)-2=\frac{k-i}{k}\ell'+2i -2 . 	$$
	This quantity is positive since $1\leq i\leq k$ and $k\geq 2$.
\end{pfcite}

In preparation to insert the vertices, we create a set of special vertex-disjoint copies of $\Khk$.

\begin{lem}\label{lem:redsets}
	There exist disjoint sets $X_i(j)\subset\tilde{V}_i^{(1)}$, $i\in [k]$, $j\in [\ell']-\{1\}$, such that for every $i\in [k]$, $j\in [\ell']-\{1\}$:
	\begin{enumerate}[(1)]
		\item $|X_i(j)|=3h\zeta n$.
		\item For every $v\in X_i(j)$, there exist two vertex-disjoint copies of $\Khk$, call them $\mathcal{R}(v)$ and $\mathcal{B}(v)$, such that
		\begin{enumerate}[(i)]
			\item $\mathcal{R}(v)$ contains $v$,
			\item $\mathcal{R}(v)$ contains $h-1$ vertices from $\tilde{V}_i^{(j)}$ and $\mathcal{B}(v)$ contains $h$ vertices from $\tilde{V}_i^{(j)}$, and
			\item for every $i'\neq i$, there exists a $j'\not\in\{1,j\}$ such that both $\mathcal{R}(v)$ and $\mathcal{B}(v)$ each have $h$ vertices from $\tilde{V}_{i'}^{(j')}$.
		\end{enumerate}
		\item The $2|X_i(j)|$ copies of $\Khk$, namely $\mathcal{R}(v)$ and $\mathcal{B}(v)$ for all $v\in X_i(j)$, are all pairwise-disjoint.
   \end{enumerate}
\end{lem}

\begin{pfcite}{Lemma~\ref{lem:redsets}}
	The proof will proceed as follows: We will have some arbitrary order on the pairs $\left\{(i,j) : i\in [k], j\in [\ell']-\{1\}\right\}$ and dynamically define
	$$ 	X=\bigcup_{(i',j')\prec (i,j)}\;\bigcup_{v\in X_{i'}(j')} \left(V(\mathcal{R}(v))\cup V(\mathcal{B}(v))\right) . 	$$
	That is, $X$ is the set of all vertices belonging to a $\mathcal{R}(v)$ or a $\mathcal{B}(v)$ for all $(i',j')$ that precede the current $(i,j)$.

	We will show that, for all $v\in X_{i}^{(j)}$ the vertex-disjoint $\mathcal{R}(v)$ and $\mathcal{B}(v)$ can be found among vertices not in $X$, as long as $|X|\leq\zeta^{1/2}L'$.

	Fix $i\in [k]$ and $j\in [\ell']$. Let $(T_1,T_2)$ be a pair of $K_k$-s in $G_r'$ from Lemma~\ref{lem:reach} for these values of $i$ and $j$.  Consider the subgraph $F$ of $G'$ induced on the sub-clusters $\tilde{V}_{i'}^{(j')}$ such that $u_{i'}^{(j')}$ form $V(T_2)$. Since $T_2$ is a $K_k$ in the reduced graph $G_r'$, every pair of sub-clusters in this subgraph is $\epsilon'$-regular with density at least $d'$.  Since $|X|\leq\zeta^{1/2}L'$, $|\tilde{V}_{i'}^{(j')}-X|\geq\frac{1}{2}|\tilde{V}_{i'}^{(j')}|$ and it follows from the definition of regularity that each pair of sub-clusters of $F-X$ is $2\epsilon'$-regular with density at least $d'-\epsilon'$. By the Key Lemma (Lemma 2.1 from~\cite{KS}), $F-X$ contains at least $3h\zeta n$ vertex-disjoint copies of $\Khk$ as long as $3h\zeta n\ll\epsilon'L'$. This is satisfied because 
	$$ 	3h\zeta n\stackrel{{\rm (P\ref{it:P2prime}}'{\rm )}}{\leq} 3h\zeta\frac{\ell'L'}{1-2\epsilon}\leq 4h\zeta\ell'L'=4h\zeta\ell DL'\leq 4h\zeta M D L'\stackrel{(\ref{eq:epsprime})}{\ll}\epsilon'L' . 	$$
	In the above inequality, we use the fact that $\ell'=\ell\cdot D\leq M\cdot D$ and $M$ and $D$ depend only on $k$ and $\epsilon$. In addition, $\zeta\ll\epsilon$.  We refer to these $3h\zeta n$ copies of $\Khk$ as \textdef{blue copies of $\Khk$} and we add their vertices to $X$.

	In a similar fashion, let $F$ now be the graph induced on the sub-clusters $\tilde{V}_{i'}^{(j')}$ such that $u_{i'}^{(j')}$ is in $V(T_1)\cup V(T_2)$.  The graph $F-X$ also satisfies the assumptions of the Key Lemma and therefore we can find $3h\zeta n$ copies of $\Khk$ in such a way that each copy has one vertex in $\tilde{V}_i^{(1)}-X$ and $h-1$ vertices in $\tilde{V}_i^{(j)}-X$.  The remaining vertices of $\Khk$ are in the sub-clusters of $V(T_1)\cap V(T_2)$. We refer to these $3h\zeta n$ copies of $\Khk$ as \textdef{red copies of $\Khk$} and add their vertices to $X$. For each red copy of $\Khk$, we put its unique vertex in $V_i^{(1)}$ into $X_i(j)$ and call this copy $\mathcal{R}(v)$. For each $v\in X_i{(j)}$, let $\mathcal{B}(v)$ be a distinct blue copy of $\Khk$ as found above.

	For this process to work, we need to ensure that $|X|\leq \zeta^{1/2}L'$ at each step.  This is true because each member of each $X_i(j)$ corresponds to two $\Khk$-s which have a total of $2hk$ vertices and, hence, $|X|\leq 2hk\sum_i\sum_j|X_i(j)|=2hk(k\ell')\cdot (3h\zeta n)=6h^2k^2\ell'\zeta n\ll\zeta^{1/2}L'$.
\end{pfcite}

We color the vertices of each $\mathcal{R}(v)$ red and the vertices of each $\mathcal{B}(v)$ blue.~\\

\subsection{Nearly-equalizing the sizes of the sub-clusters}
\label{sec:equalizing}

Let us summarize where we are: We have a designated first column (we call the first column the \textit{receptacle column} and its sub-clusters \textit{receptacle sub-clusters}) with each sub-cluster of size $L'$ and each sub-cluster having the same number of red vertices, which is at most $\zeta^{1/2}L'$.  Each such red vertex is in a different vertex-disjoint red copy of $\Khk$.  In the remaining columns, each sub-cluster has $L'$ original vertices, of which at most $\zeta^{1/2}L'$ are colored red and at most $\zeta^{1/2}L'$ are colored blue. The total number of red vertices in each $V_i$ is the same multiple of $h$. Moreover, in every column, every pair of sub-clusters is $(2\epsilon',d'/2)$-super regular. Finally, for each $i\in [k]$, there is a leftover set $\tilde{V}_i^{(0)}$ of size at most $k\epsilon'n$.

We shall now re-distribute the vertices from leftover sets $\tilde{V}_i^{(0)}$, $i\in [k]$, to non-receptacle sub-clusters in such a way that the size of leftover sets becomes $O(n)$  and each non-receptacle sub-cluster will contain exactly $h\left\lceil\left(1-d'/4\right)(L'/h)\right\rceil$ non-red vertices. These two properties will be essential for our procedure for finding perfect $\Khk$-tiling to work.~\\

We say that a vertex $v\in V_i$ \textdef{belongs in the sub-cluster $\tilde{V}_i^{(j)}$} if $v$ is adjacent to at least $(d'/2)L'$ vertices in each of the other sub-clusters $\tilde{V}_{i'}^{(j)}$, $i\neq i'$, in the $j$-th column. 

\begin{fact}\label{fact:belong}
	For every $i\in [k]$, we can partition the leftover set $\tilde{V}_i^{(0)}$ into subsets $Y_i^{(2)},\ldots,Y_i^{(\ell')}$ where, for every $j\in\{2,\ldots,\ell'\}$, the members of $Y_i^{(j)}$ belong in sub-cluster $\tilde{V}_i^{(j)}$ and 
	$$ 	\left|Y_i^{(j)}\right|\leq\frac{k\epsilon'n}{(1/k+\gamma/2)\ell'}\leq k^2\epsilon'L' . 	$$
\end{fact}

The number of red vertices in each sub-cluster may vary, but it is always less than $\zeta^{1/2}L'$.  Hence, after applying Fact~\ref{fact:belong}, the number of non-red vertices in each sub-cluster is in the interval $\left((1-\zeta^{1/2})L',(1+k^2\epsilon')L'\right)$. Fact~\ref{fact:belong} is proved in Section~\ref{sec:facts}.

Next, we wish to remove copies of $\Khk$ in such a way that the number of non-red vertices in each non-receptacle sub-cluster is the same and there are new leftover sets of size $O(\zeta n)$. This is accomplished via Lemma~\ref{lem:hatsets}. After we insert vertices via Fact~\ref{fact:belong} and remove some to create a (much smaller) leftover set via Lemma~\ref{lem:hatsets}, the sets $\tilde{V}^{(j)}_i$ will be slightly changed into sets $\hat{V}^{(j)}_i$ for $i\in [k]$ and $j\in\{0,1,\ldots,\ell'\}$.

\begin{lem}\label{lem:hatsets}
	For each $i\in [k]$, there exist disjoint vertex sets $\hat{V}_i^{(0)},\hat{V}_i^{(1)},\ldots,\hat{V}_i^{(\ell')}$ in $V_i$ such that the following occurs:
	\begin{itemize}
		\item $|\hat{V}_i^{(0)}|\leq 3h\zeta n$,
		\item $\hat{V}_i^{(1)}=\tilde{V}_i^{(1)}$, has exactly $(\ell'-1)3h\zeta n$ red vertices and exactly $L'$ vertices total,
		\item for $j\in\{2,\ldots,\ell'\}$, $\hat{V}_i^{(j)}\subset\tilde{V}_i^{(j)}$ and $\hat{V}_i^{(j)}$ contains all red and blue vertices of $\tilde{V}_i^{(j)}$,
		\item for $j\in\{2,\ldots,\ell'\}$, $\hat{V}_i^{(j)}$ contains exactly $h\left\lceil\left(1-\frac{d'}{4}\right)\frac{L'}{h}\right\rceil$ non-red vertices, and
		\item the graph induced by $V(G')-\bigcup_{i=1}^k\bigcup_{j=0}^{\ell'}\hat{V}_i^{(j)}$ is spanned by the union of vertex-disjoint copies of $\Khk$.
	\end{itemize}
\end{lem}

\begin{pfcite}{Lemma~\ref{lem:hatsets}}
	In this proof, we will remove some copies of $\Khk$ to thin the graph so that the sub-clusters satisfy the conditions above.  We shall do this by taking the reduced graph $G_r$ and creating an auxiliary graph $A_r$ and then we apply Theorem~\ref{thm:frac} to $A_r$.  From the resulting fractional $K_k$-tiling in $A_r$, we will produce a family of vertex-disjoint $\Khk$-s in $G$ that we shall remove.

	From Section~\ref{sec:partitioning}, recall that $D=D(k,\epsilon)$ was the least common multiple of a common denominator of a rational-valued solution to linear program (\ref{eq:maxLP}) over all balanced $k$-partite graphs with at most $k\cdot M=k\cdot M(k,\epsilon)$ vertices.  In a similar way, we may define $D_0=D_0(k,\epsilon,\zeta)$ to be the least common multiple of the common denominator of a rational-valued solution to linear program (\ref{eq:maxLP}) over all balanced $k$-partite graphs with at most $\frac{d'}{3h\zeta}\ell'\leq\frac{1}{3h\zeta}D(k,\epsilon)M(k,\epsilon)$ vertices in each class.

	Now we will define the auxiliary reduced graph $A_r$ by blowing up the vertices and edges of the subgraph of $G_r'$ induced by $V(G_r')-\{u_1^{(1)},u_2^{(1)},\ldots,u_k^{(1)}\}$. The number of copies of each vertex, however, will not be the same.  For $i\in [k]$ and $j\in\{2,\ldots,\ell'\}$, define $\nu(u_i^{(j)})$ to be the number of non-red vertices in sub-cluster~$\tilde{V}_i^{(j)}$.

	For $V(A_r)$, replace each vertex $u_i^{(j)}$ with the following number of copies: either the ceiling or the floor of
	$$ 	\frac{\nu(u_i^{(j)})-\lceil (1-d'/4)L'\rceil}{hD_0\lceil\zeta L'/D_0\rceil}-1 . 	$$
	The choice of ceiling or floor is made arbitrarily, but only to ensure that the resulting graph is balanced.  This is always possible because $\sum_{j=2}^{\ell'}\nu(u_i^{(j)})$ is the same for all $i\in [k]$. For $E(A_r)$, we replace each edge in $G_r'$ by a complete bipartite graph and each nonedge by an empty bipartite graph.

	First, we need to check that the number of vertices of $A_r$ is not too large. Since $(1-\zeta^{1/2})L'\leq\nu(u_i^{(j)})\leq (1+k^2\epsilon')L'$, the number of vertices in each partite set of $A_r$ is at most
	\begin{align}
		&\sum_{j=2}^{\ell'}\left\lceil\frac{\nu(u_i^{(j)})-\lceil (1-d'/4)L'\rceil}{hD_0\lceil\zeta L'/D_0\rceil}-1\right\rceil \nonumber \\
		&\qquad \leq (\ell'-1)\left\lceil\frac{(1+k^2\epsilon')L'-(1-d'/4)L'}{h\zeta L'}-1\right\rceil \nonumber \\
		&\qquad < (\ell'-1)\frac{d'/4+k^2\epsilon'}{h\zeta} . \label{it:maxnumvert}
	\end{align}
	This quantity is at most $\frac{d'}{3h\zeta}\ell'$ because $\epsilon'\ll d'$.

	Second, we need to check that each vertex $A_r$ has sufficiently large degrees in order to apply Theorem~\ref{thm:frac}. We observe that if $u$ were adjacent to $u_i^{(j)}$ in $G_r'$, then every copy of $u$ in $V(A_r)$ is adjacent to at least
	$$ 	\left\lfloor\frac{(1-\zeta^{1/2})L'-\lceil (1-d'/4)L'\rceil}{hD_0\lceil\zeta L'/D_0\rceil} -1\right\rfloor \geq\frac{d'/4-2\zeta^{1/2}}{h\zeta} 	$$
	copies of $u_i^{(j)}$ in $V(A_r)$.  So, each vertex in $V(A_r)$ is adjacent to at least
	$$ 	\left[\left(\frac{k-1}{k}+\frac{\gamma}{2}\right)\ell'-1\right]\frac{d'/4-2\zeta^{1/2}}{h\zeta}
		\geq (\ell'-1)\left(\frac{k-1}{k}+\frac{\gamma}{3}\right)\frac{d'}{4h\zeta} 	$$
	vertices in each of the other partite sets of $V(A_r)$.
	By (\ref{it:maxnumvert}), every partite set of $V(A_r)$ has size at most $(\ell'-1)\frac{d'/4+k^2\epsilon'}{h\zeta}$.  Using (\ref{eq:epsprime}), the proportion of neighbors of a vertex in $V(A_r)$ in any other vertex class is at least
	$$ 	\frac{(\ell'-1)\left(\frac{k-1}{k}+\frac{\gamma}{3}\right)\frac{d'}{4h\zeta}}{(\ell'-1)\frac{d'/4+k^2\epsilon'}{h\zeta}}\geq \frac{k-1}{k} . 	$$

	So, we can apply Theorem~\ref{thm:frac} to the auxiliary reduced graph $A_r$ and obtain an optimal solution to linear program (\ref{eq:maxLP}) with the property that $D_0 w(T)$ is an integer for every $T\in\mathcal{T}_k(A_r)$.

	As in Fact~\ref{fact:grtiling}, this implies that the blow-up graph $A_r(D_0)$ must have a perfect $K_k$-tiling.  For each $K_k$ in this tiling, we will remove $\lceil\zeta L'/D_0\rceil$ vertex-disjoint copies of $\Khk$ from the uncolored vertices of the corresponding sub-clusters of $G_r'$.
  
	It is easy to find such vertex-disjoint copies of $\Khk$ in a $k$-tuple.  Observe that every sub-cluster has at most $k^2\epsilon'L'$ uncolored vertices added to the sub-cluster.  Moreover, a set of $h D_0\lceil\zeta L'/D_0\rceil$ vertices will be removed from a sub-cluster at most $d'/(2\zeta)$ times as long as $\epsilon'\ll d'$.  So, there will always be at least $\lceil (1-d'/4)L'\rceil-k^2\epsilon'L'-h d'L'/2\geq (1-d')L'$ uncolored vertices from the original sub-cluster.  Using the Slicing Lemma (Fact~\ref{fact:slicing}), any pair of them form a $2(2\epsilon')$-regular pair.  As long as $\zeta\ll\epsilon'\ll d'$, we could apply, say, the Key Lemma from~\cite{KS} to ensure the existence of at most $\lceil\zeta L'/D_0\rceil$ vertex-disjoint copies of $\Khk$ in the $k$-tuple.
  
	So, the total number of vertices removed from sub-cluster $V_i^{(j)}$ is
	$$ 	h D_0\lceil\zeta L'/D_0\rceil\times\left\lfloor\frac{\nu(u_i^{(j)})-\lceil (1-d'/4)L'\rceil}{hD_0\lceil\zeta L'/D_0\rceil}-1\right\rceil , 	$$
	where $\lfloor\cdot\rceil$ is either the floor or ceiling of its argument.

	Removing these copies of $\Khk$ has the effect of making the number of uncolored vertices in each sub-cluster nearly identical, that is, within $h \zeta L'$ of each other. For $i\in [k]$, place into the new leftover set of $V_i$ at most $h D_0\lceil\zeta L'/D_0\rceil-1$ uncolored vertices from each sub-cluster to ensure that every sub-cluster retains either $\lceil (1-d'/4)L'\rceil$ or $\lceil (1-d'/4)L'\rceil+h D_0\lceil\zeta L'/D_0\rceil$ uncolored vertices, depending on whether the ceiling or floor function was chosen for rounding.  In the latter case, place an additional $h D_0\lceil\zeta L'/D_0\rceil$ uncolored vertices from the sub-cluster to the leftover set.

	Summarizing:
	\begin{itemize}
		\item We placed into each leftover set at most $2hD_0\lceil\zeta L'/D_0\rceil$ vertices from each sub-cluster, so each new leftover set $\hat{V}_i^{(0)}$ has a size of at most $\ell'\cdot 2hD_0\lceil\zeta L'/D_0\rceil\leq 3h\zeta n$.
		\item The sets $\tilde{V}_i^{(1)}$ are unchanged.
		\item For $j\in\{2,\ldots,\ell'\}$, $\hat{V}_i^{(j)}$ is formed by removing uncolored vertices from~$\tilde{V}_i^{(j)}$.
		\item For $j\in\{2,\ldots,\ell'\}$, the number of non-red vertices in $\hat{V}_i^{(j)}$ is explicitly prescribed to be $h\left\lceil\left(1-d'/4\right)(L'/h)\right\rceil$ because later we need it to be divisible by $h$.
		\item The vertices that are removed are all in vertex-disjoint copies of $\Khk$.
	\end{itemize}

\end{pfcite}~\\

\subsection{Inserting the leftover vertices and construction of perfect $\Khk$-tiling}
\label{sec:inserting}
We first insert the leftover vertices from $\bigcup_{i=1}^k\hat{V}_i^{(0)}$ to non-receptacle sub-clusters in such a way that we shall be able to find a perfect $\Khk$-tiling in every column using the Blow-up Lemma.  That is, each sub-cluster in the column will have the same number of vertices (divisible by $h$) and each pair of sub-clusters will be super-regular.

Suppose that vertex $w\in \hat{V}_i^{(0)}$ belongs in the sub-cluster $\tilde{V}_i^{(j)}$, $j\in\{2,\dots,\ell'\}$. We then take any $v\in X_i^{(j)}$ and the red and blue copies $\mathcal{R}(v), \mathcal{B}(v)$ of $\Khk$ guaranteed by Lemma~\ref{lem:redsets}. We uncolor the vertices of $\mathcal{R}(v)$, remove the vertices of $\mathcal{B}(v)$ from their respective sub-clusters and place $\mathcal{B}(v)$ aside to be included in the final tiling of $G$. We also add $w$ to the sub-cluster $\tilde{V}_i^{(j)}$ and remove $v$ from~$X_i^{(j)}$.

Each time this procedure is undertaken, the number of non-red vertices in each non-receptacle sub-cluster does not change and it is equal to $h\left\lceil\left(1-d'/4\right)(L'/h)\right\rceil$.

After doing this procedure for every vertex in the leftover sets, we remove all the remaining (unused) red copies of of $\Khk$ and place them aside to be included in the final tiling of $G$. The sub-clusters in the first (receptacle) column have the same number of non-red vertices as each other and the number of non-red vertices in each receptacle sub-cluster has the same congruency modulo $h$ as $n$ does.  That is, if we remove $n-h\lfloor n/h\rfloor$ non-red vertices from each receptacle sub-cluster, the remaining number of non-red vertices is divisible by $h$.

The non-red vertices in each receptacle sub-cluster form pairwise $(4\epsilon',d'/4)$-super-regular pairs, this follows from the Slicing Lemma (Fact~\ref{fact:slicing}) because no vertices were added to these sub-clusters. So we focus on the non-receptacle sub-clusters.

Throughout this proof, in every non-receptacle sub-cluster, at most $\epsilon'L'$ vertices were colored red and at most $\epsilon'L'$ red vertices will be uncolored (i.e., they become non-red).  In addition, the non-red vertices in any non-receptacle sub-cluster will have cardinality exactly $h\lceil (1-d'/4)L'/h\rceil$.  Recall that the original sub-clusters formed $(2\epsilon',d')$-super-regular pairs in each column.  There were at most $k^2\epsilon'L'$ new vertices added to each sub-cluster, each of which were adjacent to at least $(d'/2)L'$ vertices in each of the original sub-clusters of the column. The next lemma will imply that the non-red vertices in every non-receptacle column will form
super-regular pairs.
\begin{fact}\label{fact:supregadd}
	Let $(A,B)$ be an $(\epsilon_1,\delta_1)$-super-regular pair. Furthermore, let $A'\supset A$ and $B'\supset B$ be such that $|A'-A|\leq\epsilon_2|A|$ and $|B'-B|\leq\epsilon_2|B|$. If  
	\begin{itemize}
		\item every vertex in $A'-A$ has at least $\delta_2|B|$ neighbors in $B$ and
		\item every vertex in $B'-B$ has at least $\delta_2|A|$ neighbors in $A$, 
	\end{itemize}
	then the pair $(A',B')$ is $(\epsilon_0,\delta_0)$-super-regular, where $\delta_0=\frac{\min\{\delta_1,\delta_2\}}{(1+\epsilon_2)^2}$ and $\epsilon_0=\epsilon_1+\epsilon_2$.
\end{fact}

We apply Fact~\ref{fact:supregadd} with $\epsilon_1=2\epsilon'$, $\delta_1=d'$, $\epsilon_2=k\epsilon'$ and $\delta_2=d'/2$.  Consequently, we use $\epsilon_1+\epsilon_2\leq (k+2)\epsilon'\leq\sqrt{\epsilon'}$ and $\frac{\min\{\delta_1,\delta_2\}}{(1+\epsilon_2)^2} =\frac{d'/2}{(1+k\epsilon')^2}\geq d'/3$ to conclude that the augmented pairs in each column are $(\sqrt{\epsilon'},d'/3)$-super-regular.

Finally, to finish the tiling, apply the Blow-up Lemma to non-red vertices in each non-receptacle column (recall that the number of such vertices is the same and is divisible by $h$).  We can also apply the Blow-up Lemma to the non-red vertices in the receptacle column as well, because the sizes of those sets are divisible by $h$.
\begin{lem}[Blow-up Lemma, Koml\'os-S\'ark\"ozy-Szemer\'edi~\cite{KSSz97}]\label{lem:blow-up}
	Given a graph $R$ of order $r$ and positive parameters $\delta,\Delta$, there exists an $\epsilon_{\rm BL}>0$ such that the following holds: Let $N$ be an arbitrary positive integer, and let us replace the vertices of $R$ with pairwise disjoint $N$-sets $V_1,V_2,\ldots,V_r$ (blowing up). We construct two graphs on the same vertex-set $V=\bigcup V_i$.  The graph $R(N)$ the graph which is the blow-up of $R$ by $N$ and a sparser graph $G$ is constructed by replacing the edges of $R$ with some $(\epsilon_{\rm BL},\delta)$-super-regular pairs. If a graph $H$ with maximum degree $\Delta(H)\leq\Delta$ can be embedded into $R(N)$, then it can be embedded into $G$.
\end{lem}

Our $\Khk$-tiling consists of
\begin{itemize}
	\item[(i)] the copies of $\Khk$ that are outside of the sets~$\hat{V}_i^{(j)}$, as established in Lemma~\ref{lem:hatsets},
	\item[(ii)] the red copies of $\Khk$ that were not uncolored in the process of absorbing vertices from the leftover sets $\hat{V}_i^{(0)}$ to non-receptacle sub-clusters, and
	\item[(iii)] the copies of $\Khk$ found by applying the Blow-up Lemma to the non-red vertices in each column.
\end{itemize}
This is the tiling of $G$ with $n/h$ copies of $\Khk$.

What remains to show is that we can choose our constants to satisfy \eqref{eq:epsprime} so that all inequalities in our proof will be satisfied for sufficiently large $n$. Indeed, for given $\gamma>0$ and $h$, we let $d=\gamma/4$. We also set $R=K_k$, $r=k$, $\Delta=(k-1)h$  and $\delta= \gamma/12$ and apply Lemma \ref{lem:blow-up} to obtain $\epsilon_{\rm BL}$.  Now we define $\epsilon'=\min\{\epsilon^2_{\ref{lem:blow-up}}, d/(12k^2)\}$ and we let $\epsilon=\min\{(\epsilon')^5/16, d/4(k+2)\}$. Finally, we set $\zeta= 1/(12h^2k^2M(k,\eps)^2D(k,\eps)^2)$, where $M(k,\eps)$ comes from the Regularity Lemma (Theorem \ref{thm:SzemRegLem}) and $D(k, \eps)$ is defined in Section \ref{sec:partitioning}. This concludes the proof of Theorem~\ref{thm:main}.~\\


\section{Proofs of Facts}
\label{sec:facts}

For convenience, we restate the facts to be proven.

\noindent\textbf{Fact~\ref{fact:belong}.} \textit{For every $i\in [k]$, we can partition the leftover set $\tilde{V}_i^{(0)}$ into subsets $Y_i^{(2)},\ldots,Y_i^{(\ell')}$ where, for every $j\in\{2,\ldots,\ell'\}$, the members of $Y_i^{(j)}$ belong in sub-cluster $\tilde{V}_i^{(j)}$ and $|Y_i^{(j)}|\leq\tfrac{k\epsilon'n}{(1/k+\gamma/2)\ell'}\leq k^2\epsilon'L'$.}~\\

\begin{pfcite}{Fact~\ref{fact:belong}}
	First, we show that each vertex belongs in at least $(1/k+\gamma/2)\ell'$ sub-clusters. To see this, let $x$ be the number of sub-clusters in $V_{i'}$, $i'\neq i$ such that $v$ is adjacent to less than $(d'/2)L'$ vertices of that sub-cluster.  Then, since $n-\ell'L'\leq 2\epsilon n$,
	$$	x\frac{d'}{2}L'+(\ell'-x)L'+(n-\ell'L')\geq\left(\frac{k-1}{k}+\gamma\right)n . 	$$
	From this it is easy to derive that with $d',\epsilon'$ small enough relative to $\gamma$, it is the case that $x<(1/k-\gamma/2)\ell'$. By a simple union bound, the number of sub-clusters in which $v$ belongs is greater than $\ell'-(k-1)(1/k-\gamma/2)\ell'\geq (1/k+\gamma/2)\ell'$. Hence, there are at least $(1/k+\gamma/2)\ell'$ sub-clusters outside of the receptacle column in which $v$ belongs.

	Sequentially and arbitrarily assign $v\in\tilde{V}_i^{(0)}$ to $Y_i^{(j)}$ if both $v$ belongs in $\tilde{V}_i^{(j)}$ and $|Y_i^{(j)}|<\frac{k\epsilon'n}{(1/k+\gamma/2)\ell'}$. Since the size of $\tilde{V}_i^{(0)}$ is at most $k\epsilon'n$, we can always find a place for $v$.
\end{pfcite}~\\

\noindent\textbf{Fact~\ref{fact:super-slicing}.} \textit{Let $\epsilon'>0$ and $\epsilon'<d'/(2k+2)$. Let $(A_1,\ldots,A_k)$ be a $k$-tuple that is pairwise $\epsilon'$-regular of density at least $d'$ with $|A_1|=\cdots=|A_k|=L'$.  There exist subsets $A_i'\subset A_i$ for $i\in [k]$ such that $|A_i|=h\lceil (1-(k-1)\epsilon')L'/h\rceil$ and each pair of $(A_1',\ldots,A_k')$ is $(2\epsilon',d'-k\epsilon')$-super-regular (with density at least $d'-\epsilon'$).}~\\

\begin{pfcite}{Fact~\ref{fact:super-slicing}}
	We use the so-called Slicing Lemma~\cite[Fact 1.5]{KS}.
	\begin{fact}[Slicing Lemma~\cite{KS}]\label{fact:slicing}
		Given $\epsilon,\alpha,d$ such that $0<\epsilon<\alpha<1$ and $d,1-d\geq\max\{2\epsilon,\epsilon/\alpha\}$. Let $(A,B)$ be an $\epsilon$-regular pair with density $d$, $A'\subset A$ with $|A'|\geq\alpha |A|$ and $B'\subset B$ with $|B'|\geq\alpha |B|$. Then $(A',B')$ is $\epsilon_0$-regular with $\epsilon_0=\max\{2\epsilon,\epsilon/\alpha\}$ and density in $[d-\epsilon,d+\epsilon]$.
	\end{fact}

	It follows from the $\epsilon'$-regularity of $(A_i, A_j)$ that all but at most $\epsilon' |A_i|$ vertices of $A_i$ have at least $(d'-\epsilon')|A_j|$ neighbors in $A_j$. So, there is a set $A_i'\subset A_i$ of size $(1-(k-1)\epsilon')|A_i|$ such that each vertex of $A_i'$ has at least $(d'-\epsilon')|A_j|$ neighbors in $A_j$ for every $j\neq i$ and, consequently, at least $(d'-\epsilon')|A_j|-(k-1)\epsilon'|A_j|=(d'-k\epsilon')|A_j'|$ neighbors in $A_j'$ for every $j\neq i$. 

	Since $\epsilon'<d'/(2k+2)$ and $(1-(k-1)\epsilon')>1/2$, the Slicing Lemma (Fact~\ref{fact:slicing}) with $\alpha=1/2$ and $\epsilon'<d'/(2k+2)$ gives that each pair $(A_i', A_j')$ is $(2\epsilon', d-k\epsilon')$-super-regular.
\end{pfcite}~\\

\noindent\textbf{Fact~\ref{fact:supregadd}.} Let $(A,B)$ be an $(\epsilon_1,\delta_1)$-super-regular pair. Furthermore, let $A'\supset A$ and $B'\supset B$ be such that $|A'-A|\leq\epsilon_2|A|$ and $|B'-B|\leq\epsilon_2|B|$. If every vertex in $A'-A$ has at least $\delta_2|B|$ neighbors in $B$ and every vertex in $B'-B$ has at least $\delta_2|A|$ neighbors in $A$, then the pair $(A',B')$ is $(\epsilon_0,\delta_0)$-super-regular, where $\epsilon_0=\epsilon_1+\epsilon_2$ and $\delta_0=\frac{\min\{\delta_1,\delta_2\}}{(1+\epsilon_2)^2}$.~\\

\begin{pfcite}{Fact~\ref{fact:supregadd}}
	First we establish the minimum degree condition.  Each of the vertices in $A$ is adjacent to at least $\delta_1 |B|=\delta_1\frac{|B|}{|B'|}|B'|$ vertices in $B'$.  Each of the vertices in $A'-A$ is adjacent to at least $\delta_2|B|=\delta_2\frac{|B|}{|B'|}|B'|$ neighbors in $B'$. Similar conditions hold for vertices in $B'$.

	Since
	$$ 	\delta_0\leq\frac{\min\{\delta_1,\delta_2\}}{(1+\epsilon_2)^2} \leq\min\left\{\delta_1\frac{|A|}{|A'|},\delta_2\frac{|A|}{|A'|}, \delta_1\frac{|B|}{|B'|},\delta_2\frac{|B|}{|B'|}\right\} , 	$$
	each vertex $a\in A'$ has at least $\delta_0|B'|$ neighbors in $B'$ and each vertex $b\in B'$ has at least $\delta_0|A'|$ neighbors in $A'$.

	Now, consider any $X'\subseteq A'$ and $Y'\subseteq B'$ such that $|X'|\geq\epsilon_0|A'|$ and $|Y'|\geq\epsilon_0|B'|$.  Consider $X=X'-(A'-A)$ and $Y=Y'-(B'-B)$.  Note that
	$$ 	|X|\geq |X'|-\epsilon_2|A|\geq \epsilon_0|A'|-\epsilon_2|A|\geq\epsilon_1|A| . 	$$
	Similarly, $|Y|\geq\epsilon_1|B|$ and so $d(X,Y)\geq\delta_1$.

	Consequently,
	$$	d(X',Y')\geq d(X,Y)\frac{|X||Y|}{|X'||Y'|}\geq\frac{\delta_1}{(1+\epsilon_2)^{2}}\geq\delta_0 , 	$$
	and the pair is $(\epsilon_0,\delta_0)$-super-regular.
\end{pfcite}


\section{Concluding Remarks}
\label{sec:conc}
The common denominator $D=D(k,\epsilon)$ used in Section~\ref{sec:partitioning} can, in principle, be astronomically large, as it is the common denominator of values of rational-valued solutions for all balanced $k$-partite graphs on at most $M=M(k,\epsilon)$ vertices.  We chose this value for the convenience of the proof. Indeed, the constant $M$ is quite large itself and so $D$ is not so large, relatively speaking.

We could utilize a much smaller integer value of $D$ by choosing a $D$ such that if $w^*$ is the rational-valued solution of (\ref{eq:maxLP}), then for every $v\in V(G_r)$ and every $T\in\T_k(G_r)$ for which $V(T)\ni v$, we assign $\lfloor Dw^*(T)\rfloor$ vertices of $G_r'$ to copies of $T$.  Because $Dw^*(T)$ is not necessarily an integer, we end up with $D-\sum_{V(T)\ni v}\lfloor Dw^*(T)\rfloor$ unused vertices.  Choose $D$ large enough to ensure that this is always small ($O(\epsilon M^{k-1})$ suffices), and they can be placed in the leftover set.

We should also note that, asymptotically, Conjecture~\ref{conj:main} is stronger than the Hajnal-Szemer\'edi Theorem. That is, if $G$ is a graph on $kn$ vertices with minimum degree at least $\left(\frac{k-1}{k}+\gamma\right)kn$, then a random partition of the vertex set into $k$ equal parts gives a $k$-partite graph $\tilde{G}$ with $\bdelta_k(\tilde{G})\geq \left(\frac{k-1}{k}+\gamma\right)kn-O(\sqrt{n}\log n)$ and applying Conjecture~\ref{conj:main} would give a $K_k$-tiling in $\tilde{G}$ and, hence, $G$ itself.

\section{Acknowledgements}
\label{sec:ack}
The authors would like to thank the London School of Economics and Political Science.  Part of this research had been done as part of a visit by the first author. We express our appreciation to Hong-Jian Lai who has made the publications of Paul Catlin available online (\texttt{http://www.math.wvu.edu/\~{}hjlai/}). Thanks also to Bernhard von Stengel for helpful comments on linear programming and to Richard Mycroft for useful conversations. Special thanks to an anonymous referee whose comments improved the paper.

\bibliographystyle{plain}

\begin{thebibliography}{99}
%
%
\bibitem{AY} N. Alon and R. Yuster, $H$-factors in dense graphs. \textit{J. Combin. Theory Ser. B} \textbf{66} (1996), no.~2, 269--282.
\bibitem{BZ} A. Bush and Y. Zhao, Minimum degree thresholds for bipartite graph tiling. \textit{J. Graph Theory} \textbf{70} (2012), no.~1, 92--120.
\bibitem{C} P.A. Catlin, On the Hajnal-Szemer\'edi theorem on disjoint cliques. \textit{Utilitas Math.} \textbf{17} (1980), 163--177.
\bibitem{Chv} V. Chv\'atal, \textit{Linear Programming}. A series of Books in the Mathematical Sciences. \textit{W. H. Freeman and Company, New York}, 1983. xiii+478~pp.
\bibitem{CH} K. Corr\'adi and A. Hajnal, On the maximal number of independent circuits in a graph. \textit{Acta Math. Acad. Sci. Hungar.} \textbf{14} (1963), 423--439.
\bibitem{CsM} B. Csaba and M. Mydlarz, Approximate multipartite version of the Hajnal-Szemer\'edi theorem. \textit{J.~Combin. Theory Ser. B} \textbf{102} (2012), no.~2, 395--410.
\bibitem{CDeB} A. Czygrinow and L. DeBiasio, A note on bipartite graph tiling, \textit{SIAM J.~Discrete Math.}, \textbf{25} (2011), no.~4, 1477--1489.
\bibitem{Fis} E. Fischer, Variants of the Hajnal-Szemer\'edi theorem, \textit{J.~Graph Theory} \textbf{31} (1999), no.~4, 275--282.
\bibitem{Gass} S. I. Gass, \textit{Linear programming. Methods and applications.} Third edition, McGraw-Hill Book Co., New York-London-Toronto, Ont. 1969 xiii+358~pp.
\bibitem{HSz} A. Hajnal and E. Szemer\'edi, Proof of a conjecture of P. Erd\H{o}s, \textit{Combinatorial theory and its applications, II (Proc. Colloq. Balatonf\"ured, 1969)}, pp.~601--623. \textit{North-Holland}, Amsterdam, 1970.
\bibitem{Hal} P. Hall, On representation of subsets. \textit{J.~London Math. Soc.} \textbf{10} (1935), 26--30.
\bibitem{HS} J. Hladk\'y and M. Schacht, Note on bipartite graph tilings. \textit{SIAM J. Discrete Math.} \textbf{24} (2010), no. 2, 357--362.
%
%
\bibitem{Joh} R. Johansson, Triangle factors in a balanced blown-up triangle. \textit{Discrete Math.} \textbf{211} (2000), no.~1--3, 249--254.
\bibitem{KM} P. Keevash and R. Mycroft, A geometric theory for hypergraph matching, \textit{Mem. Amer. Math. Soc.} \textbf{233} (2015), no. 1098 vi+95 pp.
\bibitem{KM2} P. Keevash and R. Mycroft, A multipartite Hajnal-Szemer\'edi theorem, \textit{J.~Combin. Theory Ser. B} \textbf{114} (2015), 187--236.
\bibitem{KK} H.A. Kierstead and A.V. Kostochka, A short proof of the Hajnal-Szemer\'edi theorem on equitable colouring. \textit{Combin. Probab. Comput.} \textbf{17} (2008), no.~2, 265--270.
\bibitem{KKMSz} H.A. Kierstead, A.V. Kostochka and M. Mydlarz and E. Szemer\'edi, A fast algorithm for equitable coloring, \textit{Combinatorica} \textbf{30} (2010), no.~2, 217--224.
%
%
\bibitem{Kom} J. Koml\'os, Tiling {T}ur\'an theorems. \textit{Combinatorica} \textbf{20} (2000), no.~2, 203--218.
\bibitem{KSSz97} J. Koml\'os, G. S\'ark\"ozy and E. Szemer\'edi, Blow-up lemma. \textit{Combinatorica} \textbf{17} (1997), no.~1, 109--123.
\bibitem{KSSz01} J. Koml\'os, G. S\'ark\"ozy and E. Szemer\'edi, Proof of the Alon-Yuster conjecture. Combinatorics (Prague, 1998). \textit{Discrete Math.} \textbf{235} (2001), no. 1--3, 255--269.
\bibitem{KS} J. Koml\'os and M. Simonovits, Szemer\'edi's regularity lemma and its applications in graph theory, \textit{Combinatorics, Paul Erd\H{o}s is eighty, Vol. 2 (Keszthely, 1993)}, pp.~295--352, Bolyai Soc. Math. Stud., 2, \textit{J\'anos Bolyai Math. Soc., Budapest}, 1996.
\bibitem{Kon} D. K\"onig, Graphen und Matrizen, \textit{Mat. Lapok} \textbf{38} (1931), 116--119.
\bibitem{KO} D. K\"uhn and D. Osthus, The minimum degree threshold for perfect graph packings. \textit{Combinatorica} \textbf{29} (2009), no.~1, 65--107.
\bibitem{LM} A. Lo and K. Markstr\"om, A multipartite version of the Hajnal-Szemer\'edi theorem for graphs and hypergraphs, \textit{Combin. Probab. Comput.} \textbf{22} (2013), no.~1, 97--111.
\bibitem{MM} Cs. Magyar and R. Martin, Tripartite version of the Corr\'adi-Hajnal theorem, \textit{Discrete Math.} \textbf{254} (2002), no.~1--3, 289--308.
\bibitem{MSz} R. Martin and E. Szemer\'edi, Quadripartite version of the Hajnal-Szemer\'edi theorem, \textit{Discrete Math.} \textbf{308} (2008), no.~19, 4337--4360.
\bibitem{MZ1} R. Martin and Y. Zhao, Tiling tripartite graphs with $3$-colorable graphs, \textit{Electron. J. Combin.} \textbf{16} (2009), no. 1, Research Paper 109, 16pp.
\bibitem{MZ2} R. Martin and Y. Zhao, Tiling tripartite graphs with 3-colorable graphs: The extreme case. Preprint.
\bibitem{Sch} A. Schrijver, \textit{Theory of linear and integer programming}. Wiley-Interscience Series in Discrete Mathematics. John Wiley \& Sons, Ltd., Chichester, 1986. xii+471 pp.
%
%
\bibitem{Szem} E. Szemer\'edi, Regular partitions of graphs. \textit{Probl\`emes combinatoires et th\'eorie des graphes (Colloq. Internat. CNRS, Univ. Orsay, Orsay, 1976)}, pp. 399--401, Colloq. Internat. CNRS, 260, CNRS, Paris, 1978.
\bibitem{Z} Y. Zhao, Bipartite graph tiling, \textit{SIAM J. Discrete Math.} \textbf{23} (2009), no.~2, 888--900.
\end{thebibliography}

\end{document}